\documentclass[11pt]{article}
\usepackage{amsmath}
\usepackage{amssymb}

\usepackage{theorem}
\numberwithin{equation}{section}

\newtheorem{theorem}{Theorem}[section]

\newtheorem{lemma}[theorem]{Lemma}

\begin{document}
\title{Global well-posedness and scattering for the defocusing, cubic, nonlinear Schr{\"o}dinger equation when $n = 3$ via a linear-nonlinear decomposition}
\date{\today}
\author{Ben Dodson}
\maketitle

\noindent \textbf{Abstract:} In this paper, we prove global well-posedness and scattering for the defocusing, cubic nonlinear Schr{\"o}dinger equation when $n = 3$ and $u_{0} \in H^{s}(\mathbf{R}^{3})$, $s > 5/7$. To this end, we utilize a linear-nonlinear decomposition, similar to the decomposition used in \cite{R} for the wave equation.\vspace{5mm}

\section{Introduction}
In this paper we study the three-dimensional defocusing, cubic nonlinear Schr{\"o}dinger equation,

\begin{equation}\label{1.1}
\aligned
i u_{t} + \Delta u &= |u|^{2} u, \\
u(0,x) &= u_{0}(x) \in H^{s}(\mathbf{R}^{3}).
\endaligned
\end{equation}

\noindent $H^{s}(\mathbf{R}^{3})$ denotes the usual inhomogeneous Sobolev space.

\begin{theorem}\label{t0.1}
 If $u_{0} \in H^{s}(\mathbf{R}^{3})$, $s > \frac{1}{2}$, then there exists $T(\| u_{0} \|_{H^{s}(\mathbf{R}^{3})}) > 0$ such that $(\ref{1.1})$ is locally well - posed on $[0, T)$.
\end{theorem}

\noindent \emph{Proof:} See \cite{CaWe1}. $\Box$\vspace{5mm}

\noindent $(\ref{1.1})$ also has a local solution on $[0, T)$, $T(u_{0}) > 0$ when $u_{0} \in \dot{H}^{1/2}(\mathbf{R}^{3})$. In this case $T > 0$ depends on the profile of the initial data, not just its size. If $\| u_{0} \|_{\dot{H}^{1/2}(\mathbf{R}^{3})}$ is small, then $(\ref{1.1})$ is globally well-posed and scatters to a free solution.\vspace{5mm}

\noindent Theorem $\ref{t0.1}$ also implies that if $s > 1/2$ and a solution to $(\ref{1.1})$ only exists on a maximal interval $[0, T_{\ast})$, $T_{\ast} < \infty$, then

\begin{equation}\label{1.1.1}
\lim_{t \nearrow T_{\ast}} \| u(t) \|_{H^{s}(\mathbf{R}^{3})} = \infty.
\end{equation}

\noindent \textbf{Remark:} \cite{KM} proved that in $[0, T_{\ast})$ is a maximal interval of existence for $(\ref{1.1})$, $T_{\ast} < \infty$, then

\begin{equation}\label{1.1.2}
 \limsup_{t \nearrow T_{\ast}}	\| u(t) \|_{\dot{H}^{1/2}(\mathbf{R}^{3})} = +\infty.
\end{equation}

\noindent A solution to $(\ref{1.1})$ conserves both mass

\begin{equation}\label{1.2}
M(u(t)) = \int |u(t,x)|^{2} dx = M(u(0)),
\end{equation}

\noindent and energy

\begin{equation}\label{1.3}
E(u(t)) = \frac{1}{2} \int |\nabla u(t,x)|^{2} dx + \frac{1}{4} \int |u(t,x)|^{4} dx = E(u(0)).
\end{equation}

\noindent Thus $(\ref{1.1})$ is globally well-posed in the defocusing case when $s = 1$. \cite{GV} proved $(\ref{1.1})$ is scattering when $u_{0} \in H^{1}(\mathbf{R}^{3})$.\vspace{5mm}

\noindent \textbf{Remark:} This argument will not work for the focusing equation since $$E(u(t)) = \frac{1}{2} \int |\nabla u(t,x)|^{2} dx - \frac{1}{4} \int |u(t,x)|^{4} dx,$$ is not positive definite.\vspace{5mm}

\noindent It is conjectured that $(\ref{1.1})$ is globally well - posed in time for all data included in the local theory. \cite{CKSTT1}, extending the work of \cite{B2}, introduced the I-method. Let $I : H^{s}(\mathbf{R}^{3}) \rightarrow H^{1}(\mathbf{R}^{3})$ be a radially symmetric Fourier multiplier. By controlling the change of $E(Iu(t))$, which is no longer constant, \cite{CKSTT1} proved $(\ref{1.1})$ is globally well-posed for $s > 5/6$. In \cite{CKSTT2}, an interaction Morawetz estimate improved this result to $s > 4/5$.\vspace{5mm}

\noindent In this paper we prove

\begin{theorem}\label{t1.1}
$(\ref{1.1})$ is globally well-posed for $s > 5/7$. Additionally,

\begin{equation}\label{1.5}
\| u(t) \|_{H^{s}(\mathbf{R}^{3})} \leq C(\| u_{0} \|_{H^{s}(\mathbf{R}^{3})}),
\end{equation}

\noindent and the solution scatters. There exist $u_{\pm} \in H^{s}(\mathbf{R}^{3})$ such that

\begin{equation}\label{1.6}
\aligned
\lim_{t \rightarrow \infty} \| u(t) - e^{it \Delta} u_{+} \|_{H^{s}(\mathbf{R}^{3})} = 0, \\
\lim_{t \rightarrow \infty} \| u(-t) - e^{-it \Delta} u_{-} \|_{H^{s}(\mathbf{R}^{3})} = 0.
\endaligned
\end{equation}
\end{theorem}

\noindent To do this we combine the interaction Morawetz estimates of \cite{CKSTT2} with the linear - nonlinear decomposition. In a parallel vein, \cite{R} applied the I - method to the semilinear wave equation,

\begin{equation}\label{1.4}
\aligned
\partial_{tt} u - \Delta u = -u^{3}, \\
u(0,x) \in H^{s}(\mathbf{R}^{3}), \\
u_{t}(0,x) \in H^{s - 1}(\mathbf{R}^{3}).
\endaligned
\end{equation}

\noindent \cite{R} made a linear-nonlinear decomposition, which more effectively estimated the energy change for large times. In this paper, we will make a similar argument to prove theorem $\ref{t1.1}$.\vspace{5mm}

\noindent In $\S 2$, some preliminary facts from harmonic analysis will be mentioned. In $\S 3$, a local well-posedness result will be proved. In $\S 4$, a formula for the energy increment will be computed. In $\S 5$ a smoothing estimate using a bilinear estimate will be proved. In $\S 6$, the double-layer I-decomposition will be used to prove the theorem.

\section{Preliminaries}
Let $\phi(x)$ be a smooth, radial function,

\begin{equation}\label{0.1}
\phi(x) = \left\{
            \begin{array}{ll}
              1, & \hbox{$|x| \leq 1$;} \\
              0, & \hbox{$|x| > 2$.}
            \end{array}
          \right.
\end{equation}

\noindent Let

\begin{equation}\label{0.1.1}
 \aligned
\mathcal F(P_{\leq N} u) &= \hat{u}(\xi) \phi(\frac{\xi}{N}), \\
\mathcal F(P_{> N} u) &= \hat{u}(\xi)(1 - \phi(\frac{\xi}{N})).
\endaligned
\end{equation}

\noindent Then define the standard Littlewood - Paley decomposition,

\begin{equation}\label{0.3}
P_{N} f = u_{\leq 2N} - u_{\leq N}.
\end{equation}

\noindent We let $u_{< N} = P_{< N} u$, similarly for $u_{N}$ and $u_{> N}$. The Littlewood - Paley decomposition obeys the embedding

\begin{equation}\label{0.3.1}
\| u_{N} \|_{L^{p}(\mathbf{R}^{3})}, \| u_{< N} \|_{L^{p}(\mathbf{R}^{3})}, \| u_{> N} \|_{L^{p}(\mathbf{R}^{3})}	\lesssim_{p} \| u \|_{L^{p}(\mathbf{R}^{3})}
\end{equation}

\noindent for all $1 \leq p \leq \infty$.\vspace{5mm}

\noindent The $L^{p}$ norms obey the $l^{2}$ summation rule for $1 < p < \infty$,

\begin{equation}\label{0.3.1.1}
 \| u \|_{L^{p}(\mathbf{R}^{3})}^{2}	\sim_{p}	\sum_{j = -\infty}^{\infty}	\| u_{2^{j}} \|_{L^{p}(\mathbf{R}^{3})}^{2}.
\end{equation}

\noindent Additionally Bernstein's inequality holds. For $1 < p < \infty$,

\begin{equation}\label{0.3.2}
 \| P_{N} u \|_{L^{p}(\mathbf{R}^{3})}	\lesssim_{p} \frac{1}{N^{s}} \| u \|_{\dot{H}^{s,p}(\mathbf{R}^{3})},
\end{equation}

\noindent where $\dot{H}^{s,p}$ is the $p$ - based Sobolev space of order $s$.\vspace{5mm}

\noindent We make a high-low decomposition,

\begin{equation}\label{0.2}
u = P_{\leq N} u + P_{> N} u = u_{b} + u_{s}.
\end{equation}

\noindent \textbf{Remark:} Since we will also make a linear-nonlinear decomposition, to avoid any potential confusion we will not write $u_{b}$ for low frequencies (b for bass), rather than $u_{l}$, and $u_{s}$ (s for soprano) for high frequencies.\vspace{5mm}

\noindent The I-operator is a Fourier multiplier given by a smooth, decreasing, radially symmetric symbol,

\begin{equation}\label{0.4}
I_{N} : H^{s}(\mathbf{R}^{3}) \rightarrow H^{1}(\mathbf{R}^{3}),
\end{equation}

\begin{equation}\label{0.5}
(I_{N} f)(\xi) = m_{N}(\xi) \hat{f}(\xi),
\end{equation}

\begin{equation}\label{0.6}
m_{N}(\xi) = \left\{
           \begin{array}{ll}
             1, & \hbox{$|\xi| \leq N$;} \\
             (\frac{N}{|\xi|})^{1 - s}, & \hbox{$|\xi| > 2N$.}
           \end{array}
         \right.
\end{equation}

\noindent For the rest of the paper, we understand that $If$ refers to the function $I_{N} f$. We have the estimates,

\begin{equation}\label{0.7}
\aligned
\| \nabla Iu \|_{L_{x}^{2}(\mathbf{R}^{3})} &\lesssim N^{1 - s} \| u \|_{H^{s}(\mathbf{R}^{3})}, \\
\| u \|_{H^{s}(\mathbf{R}^{3})} &\lesssim \| Iu \|_{H^{1}(\mathbf{R}^{3})}.
\endaligned
\end{equation}

\noindent \textbf{Remark:} If $E(Iu(t))$ was a conserved quantity then $(\ref{0.7})$ would imply $(\ref{1.1})$ is globally well - posed for all $s > 1/2$. Sadly this is not true. Instead, to prove theorem $\ref{t1.1}$ we will be content to merely estimate the change of $E(Iu(t))$. This estimate occupies $\S 4$.\vspace{5mm}

\noindent By Bernstein's inequality we have

\begin{equation}\label{0.8}
\| P_{> M} u \|_{L_{t}^{p} L_{x}^{q}(J \times \mathbf{R}^{3})} \lesssim (\frac{1}{M} + \frac{1}{N^{1 - s} M^{s}}) \| \nabla Iu \|_{L_{t}^{p} L_{x}^{q}(J \times \mathbf{R}^{3})},
\end{equation}

\noindent and

\begin{equation}\label{0.9}
\| |\nabla|^{1/2} P_{> M} u \|_{L_{t}^{p} L_{x}^{q}(J \times \mathbf{R}^{3})} \lesssim (\frac{1}{M^{1/2}} + \frac{1}{N^{1 - s} M^{s - 1/2}}) \| \nabla Iu \|_{L_{t}^{p} L_{x}^{q}(J \times \mathbf{R}^{3})}.
\end{equation}

\noindent We also have the Sobolev embedding theorem, for $1 \leq p < q \leq \infty$,

\begin{equation}\label{0.10}
 \| P_{N} u \|_{L^{q}(\mathbf{R}^{3})}	\lesssim	N^{\frac{3}{p} - \frac{3}{q}} \| P_{N} u \|_{L^{p}(\mathbf{R}^{3})}.
\end{equation}

\noindent \textbf{Strichartz Estimates:} A pair $(p,q)$ will be called an admissible pair if

\begin{equation}\label{0.10}
\frac{2}{p} = 3(\frac{1}{2} - \frac{1}{q}).
\end{equation}

\noindent We will also use the Strichartz space, 

\begin{equation}\label{0.11}
\| u \|_{S^{0}(J \times \mathbf{R}^{3})} = \sup_{(p,q) \text{ admissible }} \| u \|_{L_{t}^{p} L_{x}^{q}(J \times \mathbf{R}^{3})},
\end{equation}

\noindent as well as its dual,

\begin{equation}\label{0.12}
\| u \|_{N^{0}(J \times \mathbf{R}^{3})} = \inf_{(p',q') \text{ admissible }} \| u \|_{L_{t}^{p'} L_{x}^{q'}(J \times \mathbf{R}^{3})},
\end{equation}

\noindent where $p'$, $q'$ refers to the dual exponent. See \cite{Tao} for more details. If $u(t,x)$ solves the equation

\begin{equation}\label{0.13}
\aligned
i u_{t} + \Delta u = F(t), \\
u(0,x) = u_{0},
\endaligned
\end{equation}

\begin{equation}\label{0.14}
\| u \|_{S^{0}(J \times \mathbf{R}^{3})} \lesssim \| u_{0} \|_{L^{2}(\mathbf{R}^{3})} + \| F \|_{N^{0}(J \times \mathbf{R}^{3})}.
\end{equation}\vspace{5mm}

\noindent \textbf{Bilinear Estimate}\vspace{5mm}

\noindent We will also make use of the bilinear Strichartz estimate,

\begin{lemma}\label{l4.1}
Suppose

\begin{equation}\label{4.1}
u(t,x) = e^{it \Delta} u_{0} + \int_{0}^{t} e^{i(t - \tau) \Delta} F(\tau, x) d\tau,
\end{equation}

\noindent and

\begin{equation}\label{4.2}
v(t,x) = e^{it \Delta} v_{0} + \int_{0}^{t} e^{i(t - \tau) \Delta} G(\tau, x) d\tau,
\end{equation}

\noindent with $u_{0}, F$ supported on $N \leq |\xi| \leq 2N$ and $v_{0}, G$ supported on $M \leq |\xi| \leq 2M$, $N << M$. Then for any $\delta > 0$,

\begin{equation}\label{4.3}
\aligned
\| uv \|_{L_{t,x}^{2}(J \times \mathbf{R}^{3})} \lesssim \frac{N}{M^{1/2}} &(\| u_{0} \|_{L_{x}^{2}(\mathbf{R}^{3})} + \| F \|_{L_{t}^{1} L_{x}^{2}(J \times \mathbf{R}^{3})}) \\ \times &(\| v_{0} \|_{L_{x}^{2}(\mathbf{R}^{3})} + \| G \|_{L_{t}^{1} L_{x}^{2}(J \times \mathbf{R}^{3})}).
\endaligned
\end{equation}
\end{lemma}

\noindent \emph{Proof:} See \cite{CKSTT3} for a proof of the non - endpoint result, \cite{KilVis} in the endpoint case. $\Box$\vspace{5mm}

\noindent \textbf{Interaction Morawetz Estimate}

\begin{theorem}\label{t0.1}
If $u(t,x)$ solves $(\ref{1.1})$, then

\begin{equation}\label{0.15}
\| u \|_{L_{t,x}^{4}(J \times \mathbf{R}^{3})}^{4} \lesssim \| u \|_{L_{t}^{\infty} L_{x}^{2}(J \times \mathbf{R}^{3})}^{2} \| u \|_{L_{t}^{\infty} \dot{H}_{x}^{1/2}(J \times \mathbf{R}^{3})}^{2}.
\end{equation}
\end{theorem}

\noindent \emph{Proof:} See \cite{CKSTT2}.\vspace{5mm}

\section{Local Well-posedness}
In this section we prove local well-posedness when $\| u \|_{L_{t,x}^{4}(J \times \mathbf{R}^{3})}$ is small. To that end, we prove that the norm of $u$ is controlled by the norm of $Iu$.

\begin{lemma}\label{l2.1}
\noindent If $\| u \|_{L_{t,x}^{4}(J \times \mathbf{R}^{3})} \leq \epsilon$, and $I : H^{s}(\mathbf{R}^{3}) \rightarrow H^{1}(\mathbf{R}^{3})$, $1/2 < s < 1$, then

\begin{equation}\label{2.1}
\| u \|_{L_{t}^{6} L_{x}^{9/2}(J \times \mathbf{R}^{3})} \lesssim (\epsilon^{2/3} + \frac{1}{N^{1/2}}) (1 + \| \nabla Iu \|_{S^{0}(J \times \mathbf{R}^{3})}).
\end{equation}
\end{lemma}

\noindent \emph{Proof:} Make a Littlewood-Paley decomposition. By the Sobolev embedding

\begin{equation}\label{2.2}
\| P_{\leq N} u \|_{L_{t}^{\infty} L_{x}^{6}(J \times \mathbf{R}^{3})} \lesssim \| \nabla P_{\leq N} u \|_{L_{t}^{\infty} L_{x}^{2}(J \times \mathbf{R}^{3})} \leq \| \nabla Iu \|_{S^{0}(J \times \mathbf{R}^{3})}.
\end{equation}

\noindent Interpolating this with $\| P_{\leq N} u \|_{L_{t,x}^{4}(J \times \mathbf{R}^{3})} \leq \epsilon$, $$\| P_{\leq N} u \|_{L_{t}^{6} L_{x}^{9/2}(J \times \mathbf{R}^{3})} \lesssim \epsilon^{2/3} \| \nabla Iu \|_{S^{0}(J \times \mathbf{R}^{3})}^{1/3}.$$

\noindent This takes care of the $P_{\leq N}$ part. On the other hand, when $N_{j} \geq N$,

\begin{equation}\label{2.3}
\| P_{N_{j}} u \|_{L_{t}^{6} L_{x}^{9/2}(J \times \mathbf{R}^{3})} \lesssim N_{j}^{1/2} \| u \|_{L_{t}^{6} L_{x}^{18/7}(J \times \mathbf{R}^{3})} \lesssim \frac{1}{N^{1 - s}} \frac{1}{N_{j}^{s - 1/2}} \| \nabla Iu \|_{S^{0}(J \times \mathbf{R}^{3})}.
\end{equation}

\noindent Summing over $N_{j} \gtrsim N$ gives the bound for $P_{> N} u$. $\Box$

\begin{theorem}\label{t2.2}
Suppose J is an interval such that $$\| u \|_{L_{t,x}^{4}(J \times \mathbf{R}^{3})} \leq \epsilon,$$ and $E(Iu_{0}) \leq 1$. Then $(\ref{1.1})$ is locally well-posed on J, and

\begin{equation}\label{2.4}
\| \nabla Iu \|_{S^{0}(J \times \mathbf{R}^{3})} \lesssim 1.
\end{equation}
\end{theorem}

\noindent \emph{Proof:} A solution satisfies the Duhamel formula,

\begin{equation}\label{2.5}
Iu(t,x) = e^{it \Delta} Iu_{0} + \int_{0}^{t} e^{i(t - \tau) \Delta} I(|u|^{2} u)(\tau) d\tau.
\end{equation}

\noindent Since the symbol of $\nabla I$ is strictly increasing as $|\xi| \rightarrow \infty$, $\nabla I(|u|^{2} u)$ obeys the product rule. Therefore, by $(\ref{0.14})$,

$$\| \nabla Iu \|_{S^{0}(J \times \mathbf{R}^{3})} \lesssim \| \nabla Iu_{0} \|_{L^{2}(\mathbf{R}^{3})} + \| \nabla Iu \|_{L_{t}^{2} L_{x}^{6}(J \times \mathbf{R}^{3})} \| u \|_{L_{t}^{6} L_{x}^{9/2}(J \times \mathbf{R}^{3})}^{2}$$

$$\lesssim \| \nabla Iu_{0} \|_{L^{2}(\mathbf{R}^{3})} + (\epsilon^{4/3} + \frac{1}{N}) (\| \nabla Iu \|_{S^{0}(J \times \mathbf{R}^{3})} + \| \nabla Iu \|_{S^{0}(J \times \mathbf{R}^{3})}^{3}).$$

\noindent Applying the continuity method proves the theorem. $\Box$

\section{Energy Increment}
In this section we prove an estimate on the energy increment which is well suited to making long-term estimates on the change of the modified energy.

\begin{theorem}\label{t3.1}
If u is a solution to $(\ref{1.1})$, and $J = [a,b]$ is an interval with

\begin{equation}\label{3.1}
\| u \|_{L_{t,x}^{4}(J \times \mathbf{R}^{3})}^{4} \leq \epsilon,
\end{equation}

\noindent and $E(Iu(a)) \leq 1$, then

\begin{equation}\label{3.2}
\sup_{t_{1}, t_{2} \in J} |E(Iu(t_{1})) - E(Iu(t_{2}))| \lesssim \frac{1}{N^{1-}} \| \nabla I P_{> cN} u \|_{L_{t}^{2} L_{x}^{6}(J \times \mathbf{R}^{3})}^{2} + O(\frac{1}{N^{2-}}),
\end{equation}

\noindent where $c > 0$ is some constant.
\end{theorem}

\noindent \textbf{Remark:} The energy increment in $\cite{CKSTT1}$ and $\cite{CKSTT2}$ was $$\sup_{t_{1}, t_{2} \in J} |E(Iu(t_{1})) - E(Iu(t_{2}))| \lesssim \frac{1}{N^{1-}}.$$ $(\ref{3.2})$ does not offer any advantage whatsoever for one single interval. However, we can sum $(\ref{3.2})$ over many disjoint intervals much more effectively than the estimate in $\cite{CKSTT1}$.\vspace{5mm}

\noindent \emph{Proof:} To prove this, recall the formula for energy, $(\ref{1.3})$,

\begin{equation}\label{3.3}
E(Iu(t)) = \frac{1}{2} \int |\nabla Iu(t,x)|^{2} dx + \frac{1}{4} \int |Iu(t,x)|^{4} dx.
\end{equation}

\begin{equation}\label{3.4}
\aligned
\frac{d}{dt} E(Iu(t)) = -Re \int (I \partial_{t} u(t,x)) I(|u(t,x)|^{2} \overline{u(t,x)}) dx  \\ + Re \int (I \partial_{t} u(t,x)) |Iu(t,x)|^{2} \overline{Iu(t,x)} dx.
\endaligned
\end{equation}

\noindent Since $$Iu_{t} = i \Delta Iu - i I(|u|^{2} u),$$ it suffices to estimate

\begin{equation}\label{3.5}
Re \int_{t_{1}}^{t_{2}} \int (i \Delta Iu(t,x)) [I(|u(t,x)|^{2} \overline{u(t,x)}) - |Iu(t,x)|^{2} Iu(t,x)] dx dt,
\end{equation}

\noindent and

\begin{equation}\label{3.6}
Re \int_{t_{1}}^{t_{2}} \int (i I(|u(t,x)|^{2} u(t,x)))[I(|u(t,x)|^{2} \overline{u(t,x)}) - |Iu(t,x)|^{2} \overline{Iu(t,x)}] dx dt
\end{equation}

\noindent separately.\vspace{5mm}

\noindent \textbf{The term $(\ref{3.5})$:}

\begin{equation}\label{3.7}
\aligned
(\ref{3.5}) = Re \int_{t_{1}}^{t_{2}} \int_{\Sigma} (i |\xi_{1}|^{2} \widehat{Iu}(t,\xi_{1})) [\frac{m(\xi_{2} + \xi_{3} + \xi_{4})}{m(\xi_{2}) m(\xi_{3}) m(\xi_{4})} - 1] \\
 \times \widehat{\overline{Iu}}(t,\xi_{2}) \widehat{Iu}(t,\xi_{3}) \widehat{\overline{Iu}}(t,\xi_{4}) d\xi dt,
\endaligned
\end{equation}

\noindent where $\Sigma = \{ \xi_{1} + \xi_{2} + \xi_{3} + \xi_{4} = 0 \}$ and $d \xi$ is the Lebesgue measure on the hyperplane $\Sigma$. Make a Littlewood-Paley decomposition. Without loss of generality let $N_{2} \geq N_{3} \geq N_{4}$. Consider a number of cases separately.\vspace{5mm}

\noindent \textbf{Case 1, $N_{2} << N$:} In this case $$\frac{m(\xi_{2} + \xi_{3} + \xi_{4})}{m(\xi_{2}) m(\xi_{3}) m(\xi_{4})} - 1 \equiv 0.$$

\noindent \textbf{Case 2, $N_{2} \gtrsim N$, $N_{3} << N$:}\vspace{5mm}

\noindent \emph{Case 2(a): $N_{4} \geq \frac{1}{N^{2}}$} In this case, apply the fundamental theorem of calculus. $$|\frac{m(N_{2} + N_{3} + N_{4})}{m(N_{2})} - 1| \lesssim \frac{|\nabla m(N_{2})|}{m(N_{2})} N_{3} \lesssim \frac{N_{3}}{N_{2}}.$$ Therefore,

$$(\ref{3.5}) \lesssim \sum_{N \lesssim N_{1} \sim N_{2}} \frac{N_{1}}{N_{2}^{2}} \| P_{N_{1}} \nabla Iu \|_{L_{t}^{2} L_{x}^{6}(J \times \mathbf{R}^{3})}  \| P_{N_{2}} \nabla Iu \|_{L_{t}^{2} L_{x}^{6}(J \times \mathbf{R}^{3})}$$

$$\times \sum_{\frac{1}{N^{2}} \leq N_{4} \leq N_{3} << N} N_{3} \| P_{N_{3}} Iu \|_{L_{t}^{\infty} L_{x}^{2}(J \times \mathbf{R}^{3})} \| P_{N_{4}} Iu \|_{L_{t}^{\infty} L_{x}^{6}(J \times \mathbf{R}^{3})},$$

$$\lesssim \| \nabla Iu \|_{S^{0}(J \times \mathbf{R}^{3})}^{2} \sum_{N \lesssim N_{1} \sim N_{2}} \frac{\ln(N)}{N_{1}} \|  P_{N_{1}} \nabla Iu \|_{L_{t}^{2} L_{x}^{6}(J \times \mathbf{R}^{3})}  \| P_{N_{2}} \nabla Iu \|_{L_{t}^{2} L_{x}^{6}(J \times \mathbf{R}^{3})},$$

$$\lesssim \frac{1}{N^{1-}} \| P_{> cN} \nabla Iu \|_{L_{t}^{2} L_{x}^{6}(J \times \mathbf{R}^{3})}^{2} \| \nabla Iu \|_{S^{0}(J \times \mathbf{R}^{3})}^{2}.$$

\noindent The last estimate follows by Cauchy-Schwartz and $(\ref{0.3.1.1})$.\vspace{5mm}

\noindent \emph{Case 2(b), $N_{4} \leq \frac{1}{N^{2}}$:} In this case, combine the Sobolev embedding theorem with $(\ref{3.1})$,

\begin{equation}\label{3.8}
\| P_{N_{4}} u \|_{L_{t}^{4} L_{x}^{\infty}(J \times \mathbf{R}^{3})} \lesssim N_{4}^{3/4} \| P_{N_{4}} u \|_{L_{t,x}^{4}(J \times \mathbf{R}^{3})} \lesssim \epsilon N_{4}^{3/4}.
\end{equation}

\noindent Then

$$(\ref{3.5}) \lesssim \sum_{N \lesssim N_{1} \sim N_{2}} \frac{N_{1}}{N_{2}^{2}} \| P_{N_{1}} \nabla Iu \|_{L_{t}^{4} L_{x}^{3}(J \times \mathbf{R}^{3})} \| P_{N_{2}} \nabla Iu \|_{L_{t}^{4} L_{x}^{3}(J \times \mathbf{R}^{3})}$$

$$\times \sum_{N_{4} \leq \frac{1}{N^{2}}; N_{4} \leq N_{3} << N} N_{3} \| P_{N_{3}} Iu \|_{L_{t}^{4} L_{x}^{3}(J \times \mathbf{R}^{3})} \| P_{N_{4}} Iu \|_{L_{t}^{4} L_{x}^{\infty}(J \times \mathbf{R}^{3})}$$

$$\lesssim \frac{\epsilon}{N^{5/2-}} \| \nabla Iu \|_{S^{0}(J \times \mathbf{R}^{3})}^{3}.$$\vspace{5mm}

\noindent \textbf{Case 3, $N_{2} \gtrsim N$, $N_{3} \gtrsim N$, $N_{2} \sim N_{1}$:}\vspace{5mm}

\noindent \emph{Case 3(a), $N_{4} \geq \frac{1}{N^{2}}$:} In this case make the crude estimate $$|\frac{m(N_{2} + N_{3} + N_{4})}{m(N_{2}) m(N_{3}) m(N_{4})} - 1| \lesssim \frac{1}{m(N_{3}) m(N_{4})}.$$

$$(\ref{3.5}) \lesssim \sum_{N_{1} \sim N_{2}} \frac{N_{1}}{N_{2}} \| P_{N_{1}} \nabla Iu \|_{L_{t}^{2} L_{x}^{6}(J \times \mathbf{R}^{3})}  \| P_{N_{2}} \nabla Iu \|_{L_{t}^{2} L_{x}^{6}(J \times \mathbf{R}^{3})}$$

$$\times \sum_{N_{3} \gtrsim N; N_{4} \geq \frac{1}{N^{2}}} \frac{1}{N_{3} m(N_{3}) m(N_{4})} \| P_{N_{3}} \nabla Iu \|_{L_{t}^{\infty} L_{x}^{2}(J \times \mathbf{R}^{3})} \| P_{N_{4}} Iu \|_{L_{t}^{\infty} L_{x}^{6}(J \times \mathbf{R}^{3})},$$

$$\sum_{\frac{1}{N^{2}} \leq N_{4} \leq N_{3}; N \lesssim N_{3} \leq N_{2}} \frac{1}{N_{3} m(N_{3}) m(N_{4})} \lesssim \sum_{N \lesssim N_{3} \leq N_{2}} \frac{1}{N_{3} m(N_{3})} (\ln(N) + \frac{N_{3}^{1 - s}}{N^{1 - s}}) \lesssim \frac{1}{N^{1-}}.$$

\noindent Again summing $N_{1} \sim N_{2}$ by Cauchy - Schwartz,

$$(\ref{3.5}) \lesssim \frac{1}{N^{1-}} \| P_{> cN} \nabla Iu \|_{L_{t}^{2} L_{x}^{6}(J \times \mathbf{R}^{3})}^{2} \| \nabla Iu \|_{S^{0}(J \times \mathbf{R}^{3})}^{2}.$$\vspace{5mm}

\noindent \emph{Case 3(b), $N_{4} \leq \frac{1}{N^{2}}$:} Here, $$|\frac{m(N_{2} + N_{3} + N_{4})}{m(N_{2}) m(N_{3}) m(N_{4})} - 1| \lesssim \frac{1}{m(N_{3})}.$$

\noindent Once again, use the Sobolev embedding theorem combined with $\| u \|_{L_{t,x}^{4}(J \times \mathbf{R}^{3})} \leq \epsilon$.

$$(\ref{3.5}) \lesssim \sum_{N \lesssim N_{1} \sim N_{2}} \frac{N_{1}}{N_{2}} \| P_{N_{1}} \nabla Iu \|_{L_{t}^{4} L_{x}^{3}(J \times \mathbf{R}^{3})} \| P_{N_{2}} \nabla Iu \|_{L_{t}^{4} L_{x}^{3}(J \times \mathbf{R}^{3})}$$

$$\times \sum_{N_{4} \leq \frac{1}{N^{2}}; N_{3} \gtrsim N} \frac{1}{N_{3} m(N_{3})} \| P_{N_{3}} \nabla Iu \|_{L_{t}^{4} L_{x}^{3}(J \times \mathbf{R}^{3})} \| P_{N_{4}} Iu \|_{L_{t}^{4} L_{x}^{\infty}(J \times \mathbf{R}^{3})}$$

$$\lesssim \frac{\epsilon}{N^{5/2-}} \| \nabla Iu \|_{S^{0}(J \times \mathbf{R}^{3})}^{3}.$$

\noindent \textbf{Case 4, $N_{2} \gtrsim N$, $N_{2} \sim N_{3}$, $N_{1} \lesssim N_{2}$:} \vspace{5mm}

\noindent In this case $$|\frac{m(\xi_{2} + \xi_{3} + \xi_{4})}{m(\xi_{2}) m(\xi_{3}) m(\xi_{4})} - 1| \lesssim \frac{1}{m(\xi_{2}) m(\xi_{3}) m(\xi_{4})}.$$

\noindent \emph{Case 4(a), $N_{4} \geq \frac{1}{N^{2}}$:}

\begin{equation}\label{3.8.1}
 \aligned
&(\ref{3.5}) \lesssim \sum_{N \lesssim N_{2} \sim N_{3}} \frac{1}{m(N_{2}) m(N_{3}) N_{3} N_{2}} \| P_{N_{2}} \nabla Iu \|_{L_{t}^{2} L_{x}^{6}(J \times \mathbf{R}^{3})}  \| P_{N_{3}} \nabla Iu \|_{L_{t}^{2} L_{x}^{6}(J \times \mathbf{R}^{3})}	\\
&\times \sum_{N_{1} \lesssim N_{2}; \frac{1}{N^{2}} \leq N_{4} \leq N_{3}} \frac{N_{1}}{m(N_{4})} \| P_{N_{1}} \nabla Iu \|_{L_{t}^{\infty} L_{x}^{2}(J \times \mathbf{R}^{3})} \| P_{N_{4}} Iu \|_{L_{t}^{\infty} L_{x}^{6}(J \times \mathbf{R}^{3})}.
\endaligned
\end{equation}

$$\sum_{\frac{1}{N^{2}} \leq N_{4} \leq N_{3}}	\frac{1}{m(N_{4})} \| P_{N_{4}} Iu \|_{L_{t}^{\infty} L_{x}^{6}(I \times \mathbf{R}^{3})}	\lesssim (\ln(N) + \frac{N_{3}^{1 - s}}{N^{1 - s}}) \| \nabla Iu \|_{S^{0}(J \times \mathbf{R}^{3})}.$$

\noindent Because $s > \frac{1}{2}$,

$$\sum_{N_{1} \lesssim N_{2} \sim N_{3}}	\frac{N_{1}}{m(N_{2}) m(N_{3}) N_{3}^{s} N_{2} N^{1 - s}}	\lesssim \frac{1}{N^{1-}}.$$

\noindent Therefore,

$$(\ref{3.8.1})	\lesssim \frac{1}{N^{1-}} \| P_{> cN} \nabla Iu \|_{L_{t}^{2} L_{x}^{6}(J \times \mathbf{R}^{3})}^{2} \| \nabla Iu \|_{S^{0}(J \times \mathbf{R}^{3})}^{2}.$$

\noindent \emph{Case 4(b), $N_{4} \leq \frac{1}{N^{2}}$:} As usual, use the Sobolev embedding.

$$(\ref{3.5}) \lesssim \sum_{N \lesssim N_{2} \sim N_{3}} \frac{1}{N_{2} m(N_{3}) N_{3}} \| P_{N_{2}} \nabla Iu \|_{L_{t}^{4} L_{x}^{3}(J \times \mathbf{R}^{3})} \| P_{N_{3}} \nabla Iu \|_{L_{t}^{4} L_{x}^{3}(J \times \mathbf{R}^{3})}$$

$$\times \sum_{N_{4} \leq \frac{1}{N^{2}}; N_{1} \lesssim N_{2}} N_{1} \| P_{N_{1}} \nabla Iu \|_{L_{t}^{4} L_{x}^{3}(J \times \mathbf{R}^{3})} \| P_{N_{4}} Iu \|_{L_{t}^{4} L_{x}^{\infty}(J \times \mathbf{R}^{3})}$$

$$\lesssim \frac{\epsilon}{N^{5/2-}} \| \nabla Iu \|_{S^{0}(J \times \mathbf{R}^{3})}^{3}.$$

\noindent Combining all these cases with theorem $\ref{t2.2}$ proves theorem $\ref{t3.1}$ for $(\ref{3.5})$.\vspace{5mm}

\noindent \textbf{The term $(\ref{3.6}):$} To estimate this term we use a lemma.

\begin{lemma}\label{l3.2}
\begin{equation}\label{3.9}
\| P_{M} I(|u|^{2} u) \|_{L_{t,x}^{2}(J \times \mathbf{R}^{3})} \lesssim (\frac{1}{M} + \frac{1}{N}) \| \nabla Iu \|_{S^{0}(J \times \mathbf{R}^{3})}^{3}.
\end{equation}
\end{lemma}

\noindent \emph{Proof:} Make a high-low decomposition of $u$.

\begin{equation}\label{3.9.1}
 \| \nabla I(|u_{b}|^{2} u_{b}) \|_{L_{t,x}^{2}(J \times \mathbf{R}^{3})} \lesssim \| \nabla Iu \|_{L_{t}^{2} L_{x}^{6}(J \times \mathbf{R}^{3})} \| u_{b} \|_{L_{t}^{\infty} L_{x}^{6}(J \times \mathbf{R}^{3})}^{2} \lesssim \| \nabla Iu \|_{S^{0}(J \times \mathbf{R}^{3})}^{3}.
\end{equation}

\begin{equation}\label{3.9.2}
\| \nabla I(|u_{b}|^{2} u_{s}) \|_{L_{t,x}^{2}(J \times \mathbf{R}^{3})} \lesssim \| \nabla Iu \|_{L_{t}^{2} L_{x}^{6}(J \times \mathbf{R}^{3})} \| u_{b} \|_{L_{t}^{\infty} L_{x}^{6}(J \times \mathbf{R}^{3})}^{2} \lesssim \| \nabla Iu \|_{S^{0}(J \times \mathbf{R}^{3})}^{3}.
\end{equation}

\noindent Make a similar argument for $u_{b}^{2} \bar{u}_{s}$. Next, by the Sobolev embedding theorem

$$\| I(|u_{s}|^{2} u_{b}) \|_{L_{t,x}^{2}(J \times \mathbf{R}^{3})}	\lesssim	\| \nabla I(|u_{s}|^{2} u_{b}) \|_{L_{t}^{2} L_{x}^{6/5}(J \times \mathbf{R}^{3})}$$ $$\lesssim \| \nabla Iu \|_{L_{t}^{2} L_{x}^{6}(J \times \mathbf{R}^{3})} \| u_{b} \|_{L_{t}^{\infty} L_{x}^{6}(J \times \mathbf{R}^{3})} \| u_{s} \|_{L_{t}^{\infty} L_{x}^{2}(J \times \mathbf{R}^{3})} \lesssim \frac{1}{N} \| \nabla Iu \|_{S^{0}(J \times \mathbf{R}^{3})}^{3}.$$

\noindent Make a similar argument for $u_{s}^{2} \bar{u}_{b}$. Here we applied $(\ref{0.8})$ and $(\ref{0.9})$ to show

\begin{equation}\label{3.10}
\| P_{> N} u \|_{L_{t}^{\infty} L_{x}^{2}(J \times \mathbf{R}^{3})} \lesssim \frac{1}{N}  \| \nabla Iu \|_{L_{t}^{\infty} L_{x}^{2}(J \times \mathbf{R}^{3})}.
\end{equation}

\noindent Similarly, by the Sobolev embedding, $(\ref{0.8})$, and $(\ref{0.9})$,

\begin{equation}\label{3.11}
\aligned
\| \nabla I(|u_{s}|^{2} u_{s}) \|_{L_{t}^{2} L_{x}^{6/5}(J \times \mathbf{R}^{3})} &\lesssim \| \nabla Iu \|_{L_{t}^{2} L_{x}^{6}(J \times \mathbf{R}^{3})} \| u_{s} \|_{L_{t}^{\infty} \dot{H}_{x}^{1/2}(J \times \mathbf{R}^{3})}^{2} \\ &\lesssim \frac{1}{N} \| \nabla Iu \|_{S^{0}(J \times \mathbf{R}^{3})}^{3}.
\endaligned
\end{equation}

\noindent Applying Bernstein's inequality to $(\ref{3.9.1})$ and $(\ref{3.9.2})$ proves the lemma. $\Box$\vspace{5mm}

\noindent The nonlinear term is a 6-linear term. Let $\xi_{123} = \xi_{1} + \xi_{2} + \xi_{3}$ and let $N_{123}$ be the corresponding dyadic frequency such that $N_{123} \sim |\xi_{123}|$.

\begin{equation}\label{3.13}
\aligned
(\ref{3.6}) = -\int_{t_{1}}^{t_{2}} \int_{\Sigma} i \widehat{I(|u|^{2} u)}(t,\xi_{123}) [\frac{m(\xi_{4} + \xi_{5} + \xi_{6})}{m(\xi_{4}) m(\xi_{5}) m(\xi_{6})} - 1] \\
\times \widehat{\overline{Iu}}(t,\xi_{4}) \widehat{Iu}(t,\xi_{5}) \widehat{\overline{Iu}}(t,\xi_{6}) d\xi dt,
\endaligned
\end{equation}

\noindent where $\Sigma = \{ \xi_{1} + \xi_{2} + \xi_{3} + \xi_{4} + \xi_{5} + \xi_{6} = 0 \}$ and $d\xi$ is the measure on the hyperplane. Make a Littlewood-Paley decomposition and assume without loss of generality that $N_{4} \geq N_{5} \geq N_{6}$.\vspace{5mm}

\noindent \textbf{Case 1, $N_{4} << N$:} In this case the multiplier is $\equiv 0$.\vspace{5mm}

\noindent \textbf{Case 2, $N_{4} \gtrsim N$, $N_{5} << N$:} In this case the fundamental theorem of calculus will again be used. Because $N_{5}, N_{6} << N_{4}$, $N_{123} \sim N_{4}$.\vspace{5mm}

\noindent \emph{Case 2(a), $N_{6} \geq \frac{1}{N^{2}}$:}

\begin{equation}\label{3.14}
|\frac{m(\xi_{4} + \xi_{5} + \xi_{6})}{m(\xi_{4}) m(\xi_{5}) m(\xi_{6})} - 1| \lesssim \frac{|\xi_{5}|}{|\xi_{4}|}.
\end{equation}

$$(\ref{3.6}) \lesssim \sum_{N \lesssim N_{4} \sim N_{123}} \frac{1}{N_{4}} \| P_{N_{123}} I(|u|^{2} u) \|_{L_{t,x}^{2}(J \times \mathbf{R}^{3})} \| P_{N_{4}} Iu \|_{L_{t}^{2} L_{x}^{6}(J \times \mathbf{R}^{3})}$$

$$\times \sum_{\frac{1}{N^{2}} \leq N_{6} \leq N_{5} << N} N_{5} \| P_{N_{5}} Iu \|_{L_{t}^{\infty} L_{x}^{6}(J \times \mathbf{R}^{3})} \| P_{N_{6}} Iu \|_{L_{t}^{\infty} L_{x}^{6}(J \times \mathbf{R}^{3})}$$

$$\lesssim \ln(N) N \| \nabla Iu \|_{S^{0}(J \times \mathbf{R}^{3})}^{6} \sum_{N \lesssim N_{4} \sim N_{123}} \frac{1}{N_{4}^{2}} (\frac{1}{N_{123}} + \frac{1}{N}) \lesssim \frac{1}{N^{2-}} \| \nabla Iu \|_{S^{0}(J \times \mathbf{R}^{3})}^{6}.$$

\noindent \emph{Case 2(b): $N_{6} \leq \frac{1}{N^{2}}$:} As before we use the Sobolev embedding $$\| P_{N_{6}} Iu \|_{L_{t}^{4} L_{x}^{\infty}(J \times \mathbf{R}^{3})} \lesssim N_{6}^{3/4} \| P_{N_{6}} Iu \|_{L_{t,x}^{4}(J \times \mathbf{R}^{3})} \lesssim \epsilon N_{6}^{3/4}.$$

$$(\ref{3.6}) \lesssim \sum_{N \lesssim N_{123} \sim N_{4}} \frac{1}{N_{4}} \| P_{N_{123}} I(|u|^{2} u) \|_{L_{t,x}^{2}(J \times \mathbf{R}^{3})} \| P_{N_{4}} Iu \|_{L_{t}^{4} L_{x}^{3}(J \times \mathbf{R}^{3})}$$

$$\times \sum_{N_{6} \leq N_{5} << N; N_{6} \leq \frac{1}{N^{2}}} N_{5} \| P_{N_{5}} Iu \|_{L_{t}^{\infty} L_{x}^{6}(J \times \mathbf{R}^{3})} \| P_{N_{6}} Iu \|_{L_{t}^{4} L_{x}^{\infty}(J \times \mathbf{R}^{3})}$$

$$\lesssim \epsilon \sum_{N \lesssim N_{123} \sim N_{4}} (\frac{1}{N} + \frac{1}{N_{123}}) \frac{N}{N_{4}^{2}} \frac{1}{N^{3/2}} \| \nabla Iu \|_{S^{0}(J \times \mathbf{R}^{3})}^{5}$$

$$\lesssim \frac{\epsilon}{N^{7/2-}} \| \nabla Iu \|_{S^{0}(J \times \mathbf{R}^{3})}^{5}.$$\vspace{5mm}

\noindent \textbf{Case 3, $N_{5} \gtrsim N$, $N_{4} \sim N_{123}$:} Here make the crude estimate,

\begin{equation}\label{3.15}
|\frac{m(\xi_{4} + \xi_{5} + \xi_{6})}{m(\xi_{4}) m(\xi_{5}) m(\xi_{6})} - 1| \lesssim \frac{1}{m(\xi_{5}) m(\xi_{6})}.
\end{equation}

\noindent \emph{Case 3(a), $N_{6} \geq \frac{1}{N^{2}}$:}

$$(\ref{3.6}) \lesssim \sum_{N \lesssim N_{4} \sim N_{123}} \| P_{N_{123}} I(|u|^{2} u) \|_{L_{t,x}^{2}(J \times \mathbf{R}^{3})} \| P_{N_{4}} Iu \|_{L_{t}^{2} L_{x}^{6}(J \times \mathbf{R}^{3})}$$

$$\times \sum_{\frac{1}{N^{2}} \leq N_{6} \leq N_{5}; N \lesssim N_{5}} \frac{1}{m(N_{5}) m(N_{6})} \| P_{N_{5}} Iu \|_{L_{t}^{\infty} L_{x}^{6}(J \times \mathbf{R}^{3})} \| P_{N_{6}} Iu \|_{L_{t}^{\infty} L_{x}^{6}(J \times \mathbf{R}^{3})}.$$

$$\sum_{\frac{1}{N^{2}} \leq N_{6} \leq N_{5}} \frac{\| P_{N_{6}} Iu \|_{L_{t}^{\infty} L_{x}^{6}(J \times \mathbf{R}^{3})}}{m(N_{6})} \lesssim (\ln(N) + \frac{N_{5}^{1 - s}}{N^{1 - s}}) \| \nabla Iu \|_{S^{0}(J \times \mathbf{R}^{3})}.$$

$$\| \nabla Iu \|_{S^{0}(J \times \mathbf{R}^{3})} \sum_{N \lesssim N_{5} \lesssim N_{4}} \frac{1}{m(N_{5})} \| P_{N_{5}} Iu \|_{L_{t}^{\infty} L_{x}^{6}(J \times \mathbf{R}^{3})} (\ln(N) + \frac{N_{5}^{1 - s}}{N^{1 - s}}) $$ $$\lesssim (\ln(N)^{2} + \frac{N_{4}^{2(1 - s)}}{N^{2(1 - s)}}) \| \nabla Iu \|_{S^{0}(J \times \mathbf{R}^{3})}^{2}.$$

\noindent So in this case,

$$(\ref{3.6}) \lesssim \| \nabla Iu \|_{S^{0}(J \times \mathbf{R}^{3})}^{6} \sum_{N \lesssim N_{4} \sim N_{123}} \frac{1}{N} (\ln(N)^{2} + \frac{N_{4}^{2(1 - s)}}{N^{2(1 - s)}}) \frac{1}{N_{4}} \lesssim \frac{1}{N^{2-}} \| \nabla Iu \|_{S^{0}(J \times \mathbf{R}^{3})}^{6}.$$

\noindent \emph{Case 3(b), $N_{6} \leq \frac{1}{N^{2}}$:} $$|\frac{m(\xi_{4} + \xi_{5} + \xi_{6})}{m(\xi_{4}) m(\xi_{5}) m(\xi_{6})} - 1| \lesssim \frac{1}{m(\xi_{5})}.$$

$$(\ref{3.6}) \lesssim \sum_{N \lesssim N_{123} \sim N_{4}} \| P_{N_{123}} I(|u|^{2} u) \|_{L_{t,x}^{2}(J \times \mathbf{R}^{3})} \| P_{N_{4}} Iu \|_{L_{t}^{4} L_{x}^{3}(J \times \mathbf{R}^{3})}$$

$$\times \sum_{N_{5} \gtrsim N; N_{6} \leq \frac{1}{N^{2}}} \frac{1}{m(N_{5})} \| P_{N_{5}} Iu \|_{L_{t}^{\infty} L_{x}^{6}(J \times \mathbf{R}^{3})} \| P_{N_{6}} Iu \|_{L_{t}^{4} L_{x}^{\infty}(J \times \mathbf{R}^{3})}$$

$$\lesssim \sum_{N \lesssim N_{123} \sim N_{4}} \epsilon (\frac{1}{N} + \frac{1}{N_{123}}) \frac{1}{N_{4}} \frac{N_{4}^{1 - s}}{N^{1 - s}} \frac{1}{N^{3/2}} \| \nabla Iu \|_{S^{0}(J \times \mathbf{R}^{3})}^{5} $$

$$\lesssim \frac{\epsilon}{N^{7/2-}} \| \nabla Iu \|_{S^{0}(J \times \mathbf{R}^{3})}^{5}.$$\vspace{5mm}

\noindent \textbf{Case 4, $N_{5} \gtrsim N$, $N_{4} \sim N_{5}$, $N_{123} \lesssim N_{4}$:} Make the crude estimate $$|\frac{m(\xi_{4} + \xi_{5} + \xi_{6})}{m(\xi_{4}) m(\xi_{5}) m(\xi_{6})} - 1| \lesssim \frac{1}{m(\xi_{4}) m(\xi_{5}) m(\xi_{6})}.$$\vspace{5mm}

\noindent \emph{Case 4(a), $N_{6} \geq \frac{1}{N^{2}}$:}

$$(\ref{3.6}) \lesssim \sum_{N \lesssim N_{4} \sim N_{5}} \frac{1}{m(N_{4}) m(N_{5})} \| P_{N_{4}} I(u) \|_{L_{t}^{4} L_{x}^{3}(J \times \mathbf{R}^{3})} \| P_{N_{5}} Iu \|_{L_{t}^{4} L_{x}^{3}(J \times \mathbf{R}^{3})}$$

$$\times [\sum_{\frac{1}{N^{2}} \leq N_{6} \leq N_{5}; \frac{1}{N^{2}} \leq N_{123} \lesssim N_{4}} \frac{1}{ m(N_{6})} \| P_{N_{123}} I(|u|^{2} u) \|_{L_{t}^{2} L_{x}^{6} (J \times \mathbf{R}^{3})} \| P_{N_{6}} Iu \|_{L_{t}^{\infty} L_{x}^{6}(J \times \mathbf{R}^{3})}$$

$$ + \sum_{\frac{1}{N^{2}} \leq N_{6} \leq N_{5}; N_{123} \leq \frac{1}{N^{2}}} \frac{1}{ m(N_{6})} \| P_{N_{123}} I(|u|^{2} u) \|_{L_{t}^{2} L_{x}^{\infty} (J \times \mathbf{R}^{3})} \| P_{N_{6}} Iu \|_{L_{t}^{\infty} L_{x}^{3}(J \times \mathbf{R}^{3})}]$$

$$\lesssim \| \nabla Iu \|_{S^{0}(J \times \mathbf{R}^{3})}^{6} \sum_{N \lesssim N_{4} \sim N_{5}} \frac{1}{N_{4} N_{5} m(N_{4}) m(N_{5}) }$$

$$\times [\sum_{\frac{1}{N^{2}} \leq N_{6} \leq N_{5}; \frac{1}{N^{2}} \leq N_{123} \lesssim N_{4}} \frac{1}{m(N_{6})} (1 + \frac{N_{123}}{N}) + \sum_{\frac{1}{N^{2}} \leq N_{6} \leq N_{5}; N_{123} \leq \frac{1}{N^{2}}} (\frac{N_{123}^{3/2}}{N} + N_{123}^{1/2}) \frac{1}{m(N_{6}) N_{6}^{1/2}}]$$

$$\lesssim \frac{1}{N^{2-}} \| \nabla Iu \|_{S^{0}(J \times \mathbf{R}^{3})}^{6}.$$

\noindent \emph{Case 4(b), $N_{6} \leq \frac{1}{N^{2}}$:}

$$(\ref{3.6}) \lesssim \sum_{N \lesssim N_{4} \sim N_{5}} \frac{1}{m(N_{4}) m(N_{5})} \| P_{N_{4}} Iu \|_{L_{t}^{\infty} L_{x}^{2}(J \times \mathbf{R}^{3})} \| P_{N_{5}} Iu \|_{L_{t}^{4} L_{x}^{3}(J \times \mathbf{R}^{3})}$$

$$\times [\sum_{N_{123} \leq \frac{1}{N^{2}}; N_{6} \leq \frac{1}{N^{2}}}  \| P_{N_{123}} Iu \|_{L_{t}^{2} L_{x}^{\infty}(J \times \mathbf{R}^{3})} \| P_{N_{6}} Iu \|_{L_{t}^{4} L_{x}^{6}(J \times \mathbf{R}^{3})}$$

$$+ \sum_{\frac{1}{N^{2}} \leq N_{123} \lesssim N_{4}; N_{6} \leq \frac{1}{N^{2}}}  \| P_{N_{123}} Iu \|_{L_{t}^{2} L_{x}^{6}(J \times \mathbf{R}^{3})} \| P_{N_{6}} Iu \|_{L_{t}^{4} L_{x}^{\infty}(J \times \mathbf{R}^{3})}]$$

$$\lesssim \| \nabla Iu \|_{S^{0}(J \times \mathbf{R}^{3})}^{5} \sum_{N \lesssim N_{4} \sim N_{5}} \frac{1}{N_{4} N_{5} m(N_{4}) m(N_{5})} [\frac{\epsilon}{N^{3/2}} + (\ln(N) + \frac{N_{4}}{N}) \frac{\epsilon}{N^{3/2}}]$$

$$\lesssim \frac{\epsilon}{N^{7/2-}} \| \nabla Iu \|_{S^{0}(J \times \mathbf{R}^{3})}^{5}.$$

\noindent This concludes the proof of theorem $\ref{t3.1}$. $\Box$

\section{A Smoothing Estimate}
In this section we take advantage of lemma $\ref{l4.1}$ to prove a smoothing-type estimate for the Duhamel term.

\begin{lemma}\label{l4.2}
\noindent Take $N_{j} \leq N$, if 

\begin{equation}\label{4.4.0}
\| u \|_{L_{t,x}^{4}(J \times \mathbf{R}^{3})} \leq \epsilon,
\end{equation}

\noindent then

\begin{equation}\label{4.4}
\|P_{N_{j}} (|u|^{2} u) \|_{L_{t}^{1} L_{x}^{2}(J \times \mathbf{R}^{3})} \lesssim \frac{1}{N_{j}} \| P_{N_{j}} \nabla I(|u|^{2} u) \|_{L_{t}^{1} L_{x}^{2}(J \times \mathbf{R}^{3})} \lesssim \frac{1}{N_{j}} \| \nabla Iu \|_{S^{0}(J \times \mathbf{R}^{3})}^{3}.
\end{equation}
\end{lemma}

\noindent \emph{Proof:} The first inequality is Bernstein's inequality. Because $m(\xi) |\xi|$ is increasing,

$$\| P_{N_{j}} \nabla I(|u|^{2} u) \|_{L_{t}^{1} L_{x}^{2}(J \times \mathbf{R}^{3})} \lesssim \| \nabla Iu \|_{L_{t}^{2} L_{x}^{6}(J \times \mathbf{R}^{3})}$$ $$\times (\| P_{\leq 1} u \|_{L_{t}^{4} L_{x}^{6}(J \times \mathbf{R}^{3})}^{2} + \| P_{> 1} u \|_{L_{t}^{4} L_{x}^{6}(J \times \mathbf{R}^{3})}^{2}).$$ By the Sobolev embedding theorem and $(\ref{4.4.0})$,

$$\| P_{\leq 1} u \|_{L_{t}^{4} L_{x}^{6}(J \times \mathbf{R}^{3})} \lesssim \| u \|_{L_{t,x}^{4}(J \times \mathbf{R}^{3})} \leq \epsilon.$$ On the other hand,

$$\| P_{N_{k}} u \|_{L_{t}^{4} L_{x}^{6}(J \times \mathbf{R}^{3})} \lesssim N_{k}^{1/2} \| P_{N_{k}} u \|_{S^{0}(J \times \mathbf{R}^{3})}.$$

\noindent Therefore,

$$\| P_{> 1} u \|_{L_{t}^{4} L_{x}^{6}(J \times \mathbf{R}^{3})} \lesssim \sum_{1 \leq N_{k} \leq N} \frac{1}{N_{k}^{1/2}} \| \nabla Iu \|_{S^{0}(J \times \mathbf{R}^{3})}$$ $$+ \sum_{N_{k} > N} \frac{1}{N_{k}^{s - 1/2} N^{1 - s}} \| \nabla Iu \|_{S^{0}(J \times \mathbf{R}^{3})} \lesssim \| \nabla Iu \|_{S^{0}(J \times \mathbf{R}^{3})}.$$

\noindent $\Box$

\begin{theorem}\label{t4.3}
Suppose $J = [0, T]$ is an interval with 

\begin{equation}\label{4.4.1}
\| u \|_{L_{t,x}^{4}(J \times \mathbf{R}^{3})} \leq \epsilon,
\end{equation}

\noindent and $\| \nabla Iu_{0} \|_{L^{2}(\mathbf{R}^{3})} \leq 1$. The solution to $(\ref{0.1})$ on $[0, T]$ can be split into a linear piece and a nonlinear piece,

\begin{equation}\label{4.5}
u(t) = e^{it \Delta} u_{0} + \int_{0}^{t} e^{i(t - \tau) \Delta} (|u|^{2} u)(\tau) d\tau = u^{l}(t) + u^{nl}(t),
\end{equation}

\noindent with

\begin{equation}\label{4.6}
\| P_{> N} \nabla Iu^{nl} \|_{S^{0}(J \times \mathbf{R}^{3})} \lesssim \frac{1}{N^{1/2-}} (1 + \| \nabla Iu \|_{S^{0}(J \times \mathbf{R}^{3})}^{7}),
\end{equation}

\noindent and

\begin{equation}\label{4.7}
\| P_{> N} \nabla Iu^{nl} \|_{L_{t}^{\infty} L_{x}^{2}(J \times \mathbf{R}^{3})} \lesssim \frac{1}{N^{1-}} (1 + \| \nabla Iu \|_{S^{0}(J \times \mathbf{R}^{3})}^{9}).
\end{equation}
\end{theorem}

\noindent \emph{Proof:} Make a high-low decomposition of $u$, $u = u_{b} + u_{s}$ with $u_{b} = P_{\leq N/20} u$.\vspace{5mm}

\noindent Since $P_{>N} (|u_{b}|^{2} u_{b}) \equiv 0$, it suffices to consider $O(u^{2} u_{s})$. Because $|\xi| m(|\xi|)$ is increasing as $|\xi| \rightarrow \infty$,

$$\| \nabla I(|u_{b}|^{2} u_{s}) \|_{N^{0}(J \times \mathbf{R}^{3})} \lesssim \| (\nabla Iu_{s}) |u_{b}|^{2} \|_{L_{t}^{4/3} L_{x}^{3/2}(J \times \mathbf{R}^{3})}$$ $$\lesssim \| (\nabla Iu_{s}) u_{b} \|_{L_{t,x}^{2}(J \times \mathbf{R}^{3})} \| u_{b} \|_{L_{t}^{4} L_{x}^{6}(J \times \mathbf{R}^{3})}.$$

\noindent By Sobolev embedding, $(\ref{4.4.1})$, Strichartz estimates, and $\dot{H}^{1/2} \subset \dot{H}^{1}$ when $|\xi| \geq 1$,

$$\| u_{b} \|_{L_{t}^{4} L_{x}^{6}(J \times \mathbf{R}^{3})} \leq \| u_{\leq 1} \|_{L_{t}^{4} L_{x}^{6}(J \times \mathbf{R}^{3})} + \| u_{\geq 1} \|_{L_{t}^{4} L_{x}^{6}(J \times \mathbf{R}^{3})}$$

\begin{equation}\label{4.7.0}
\lesssim \epsilon + \| \nabla Iu \|_{S^{0}(J \times \mathbf{R}^{3})}.
\end{equation}

\noindent Next,

$$\| (\nabla Iu_{s}) (P_{\leq N^{-2}} u_{b}) \|_{L_{t,x}^{2}(J \times \mathbf{R}^{3})} \lesssim \| P_{\leq N^{-2}} u_{b} \|_{L_{t}^{4} L_{x}^{6}(J \times \mathbf{R}^{3})} \| \nabla Iu_{s} \|_{L_{t}^{4} L_{x}^{3}(J \times \mathbf{R}^{3})}$$

$$\lesssim \epsilon N^{-1/2} \| \nabla Iu \|_{S^{0}(J \times \mathbf{R}^{3})}.$$

\noindent Finally, estimate $$\| (\nabla Iu_{s}) (P_{> N^{-2}} u_{b}) \|_{L_{t,x}^{2}(J \times \mathbf{R}^{3})}$$ using the bilinear estimates in $(\ref{4.3})$ and lemma $\ref{l4.2}$,

$$\| (\nabla Iu_{s}) P_{> N^{-2}} u_{b} \|_{L_{t,x}^{2}(J \times \mathbf{R}^{3})} \lesssim (\sum_{N^{-2} \leq N_{k} \leq N/20} \frac{1}{N_{k}} \frac{N_{k}}{N^{1/2}}) (\| \nabla Iu \|_{S^{0}(J \times \mathbf{R}^{3})}^{2} + \| \nabla Iu \|_{S^{0}(J \times \mathbf{R}^{3})}^{6})$$

$$\lesssim \frac{1}{N^{1/2-}} (\| \nabla Iu \|_{S^{0}(J \times \mathbf{R}^{3})}^{2} + \| \nabla Iu \|_{S^{0}(J \times \mathbf{R}^{3})}^{6}).$$

\noindent Therefore,

\begin{equation}\label{4.7.1}
\| (\nabla Iu_{s}) u_{b} \|_{L_{t,x}^{2}(J \times \mathbf{R}^{3})} \lesssim \frac{1}{N^{1/2-}} (1 + \| \nabla Iu \|_{S^{0}(J \times \mathbf{R}^{3})}^{6}),
\end{equation}

\noindent which combined with $(\ref{4.7.0})$ takes care of the term $I(|u_{b}|^{2} u_{s})$. The term $I(u_{b}^{2} \bar{u}_{s})$ can be estimated in a similar manner.\vspace{5mm}
 
\noindent The other terms are easier to estimate. $$\| \nabla I(|u_{h}|^{2} u_{l}) \|_{L_{t}^{2} L_{x}^{6/5}(J \times \mathbf{R}^{3})} \lesssim \| \nabla Iu \|_{L_{t}^{2} L_{x}^{6}(J \times \mathbf{R}^{3})} \| u_{h} \|_{L_{t}^{\infty} L_{x}^{2}(J \times \mathbf{R}^{3})} \| u_{l} \|_{L_{t}^{\infty} L_{x}^{6}(J \times \mathbf{R}^{3})}$$ $$ \lesssim \frac{1}{N^{1-}} \| \nabla Iu \|_{S^{0}(J \times \mathbf{R}^{3})}^{3}.$$

\noindent A similar calculation can be made for $\bar{u}_{b} u_{s}^{2}$. Finally, $$\| \nabla I(|u_{h}|^{2} u_{h}) \|_{L_{t}^{1} L_{x}^{2}(J \times \mathbf{R}^{3})} \lesssim \| \nabla Iu \|_{L_{t}^{2} L_{x}^{6}(J \times \mathbf{R}^{3})} \| u_{h} \|_{L_{t}^{4} L_{x}^{6}(J \times \mathbf{R}^{3})}^{2} \lesssim \frac{1}{N^{1-}} \| \nabla Iu \|_{S^{0}(J \times \mathbf{R}^{3})}^{3}.$$ This finishes the proof of $(\ref{4.6})$. To prove $(\ref{4.7})$ it only remains to show

\begin{equation}\label{4.11}
\| \int_{0}^{t} e^{i(t - \tau) \Delta} P_{> N}(\nabla I(u_{b}^{2} u_{s})(\tau)) d\tau \|_{L_{t}^{\infty} L_{x}^{2}(J \times \mathbf{R}^{3})} \lesssim \frac{1}{N^{1-}} (1 + \| \nabla Iu \|_{S^{0}(J \times \mathbf{R}^{3})}^{9}).
\end{equation}

\noindent Take a function $f(t,x)$ supported on $|\xi| \geq \frac{N}{4}$ such that $$\| f(t,x) \|_{L_{t}^{1} L_{x}^{2}(J \times \mathbf{R}^{3})} = 1.$$ By duality, estimating $(\ref{4.11})$ is equivalent to estimating

\begin{equation}\label{4.11.1}\int_{J} \langle \int_{0}^{t} e^{i(t - \tau) \Delta} (\nabla I(|u_{b}|^{2} u_{s})(\tau)) d\tau, f(t,x) \rangle dt,
\end{equation}

\noindent for all such $f(t,x)$. By Fubini's theorem,

$$(\ref{4.11.1}) = \int_{J} \langle (\nabla I(|u_{b}|^{2} u_{s})(\tau)), \int_{\tau}^{T} e^{i(\tau - t) \Delta} f(t,x) dt \rangle d\tau.$$

\noindent Let

$$\int_{\tau}^{T} e^{i(\tau - t) \Delta} f(t,x) dt = v(\tau, x),$$ where $v(\tau, x)$ solves the partial differential equation

\begin{equation}\label{4.12}
\aligned
i v_{\tau} - \Delta v &= -f(\tau, x) \\
v(T) &= 0.
\endaligned
\end{equation}

$$\int_{J} \langle (\nabla I(|u_{b}|^{2} u_{s})(\tau)), v(\tau) \rangle d\tau \lesssim \| (\nabla Iu_{s}) u_{b} \|_{L_{t,x}^{2}(J \times \mathbf{R}^{3})} \| (v) u_{b} \|_{L_{t,x}^{2}(J \times \mathbf{R}^{3})}.$$ By $(\ref{4.7.1})$, $$\| (\nabla Iu_{s}) u_{b} \|_{L_{t,x}^{2}(J \times \mathbf{R}^{3})} \lesssim \frac{1}{N^{1/2-}} (1 + \| \nabla Iu \|_{S^{0}(J \times \mathbf{R}^{3})}^{6}).$$ Similarly,

$$\| v u_{b} \|_{L_{t,x}^{2}(J \times \mathbf{R}^{3})} \lesssim \frac{1}{N^{1/2-}} (1 + \| \nabla Iu \|_{S^{0}(J \times \mathbf{R}^{3})}^{3}).$$

 $\Box$

\section{Double Layer I-decomposition} Now we finally have enough tools to prove the main theorem.

\begin{theorem}\label{t5.1}
Suppose $s > 5/7$. Then $(\ref{1.1})$ is globally well-posed on $[0, \infty)$. Moreover, $\| u(t) \|_{H^{s}(\mathbf{R}^{3})} \leq C(s, \| u_{0} \|_{H^{s}(\mathbf{R}^{3})})$, and there is scattering.
\end{theorem}

\noindent \emph{Proof:} If $u(t,x)$ solves $(\ref{1.1})$ on $[0, T]$, then $\frac{1}{\lambda} u(\frac{t}{\lambda^{2}}, \frac{x}{\lambda})$ solves $(\ref{1.1})$ on $[0, \lambda^{2} T]$. This scaling leaves the $\dot{H}^{1/2}$ norm invariant. We will denote the rescaled solution $u_{\lambda}(t,x)$.

\begin{equation}\label{5.0}
\| u_{\lambda}(0,x) \|_{L^{2}(\mathbf{R}^{3})} = \lambda^{1/2} \| u_{0} \|_{L^{2}(\mathbf{R}^{3})},
\end{equation}

\begin{equation}\label{5.0.1}
\| u_{\lambda}(0,x) \|_{\dot{H}^{1}(\mathbf{R}^{3})} = \lambda^{-1/2} \| u_{0} \|_{\dot{H}^{1}(\mathbf{R}^{3})}.
\end{equation}

\noindent Combining the scaling identities with the estimates on the $I$ - operator, $(\ref{0.7})$,

$$\int |\nabla I u_{0, \lambda}(x)|^{2} dx \leq \frac{C N^{2(1 - s)}}{\lambda^{2s - 1}} \| u_{0} \|_{H^{s}(\mathbf{R}^{3})}^{2}.$$ $$\int |Iu_{0, \lambda}(x)|^{4} dx \leq \frac{C N^{3 - 4s}}{\lambda^{4s - 2}} \| u_{0} \|_{H^{s}(\mathbf{R}^{3})}^{4}.$$ Choose $\lambda \sim N^{\frac{1 - s}{s - 1/2}}$ so that $E(Iu_{0}) \leq \frac{1}{2}$. Define a set

\begin{equation}\label{5.1}
W = \{ t : E(Iu_{\lambda}(t)) \leq \frac{9}{10} \}.
\end{equation}

\noindent Since $0 \in W$, $W \neq \emptyset$. Also, by the dominated convergence theorem, $W$ is closed. So it remains to prove $W$ is open in $[0, \infty)$.\vspace{5mm}

\noindent If $W = [0, T]$, then by continuity of $E(Iu(t))$ there exists $\delta > 0$ such that $E(Iu_{\lambda}(t)) \leq 1$ on $[0, T + \delta]$.

\begin{equation}\label{5.2}
\| P_{\leq N} u \|_{L_{t}^{\infty} \dot{H}^{1/2}(J \times \mathbf{R}^{3})} \leq \| u_{0} \|_{L^{2}(\mathbf{R}^{3})}^{1/2} \| \nabla Iu \|_{L_{t}^{\infty} L_{x}^{2}(J \times \mathbf{R}^{3})}^{1/2}.
\end{equation}

\noindent Also,

\begin{equation}\label{5.3}
\| P_{> N} u \|_{L_{t}^{\infty} \dot{H}_{x}^{1/2}(J \times \mathbf{R}^{3})} \leq \frac{1}{N^{1/2}} \| \nabla Iu \|_{L_{t}^{\infty} L_{x}^{2}(J \times \mathbf{R}^{3})}.
\end{equation}

\noindent Combining the interaction Morawetz estimate $(\ref{0.15})$, $(\ref{5.2})$ and $(\ref{5.3})$,

\begin{equation}\label{5.4}
\| u_{\lambda} \|_{L_{t,x}^{4}([0, T + \delta] \times \mathbf{R}^{3})}^{4} \leq CN^{\frac{3(1 - s)}{2s - 1}}.
\end{equation}

\noindent Partition $[0, T + \delta]$ into $\sim N^{\frac{3(1 - s)}{2s - 1}}$ subintervals with $\| u_{\lambda} \|_{L_{t,x}^{4}(J_{k} \times \mathbf{R}^{3})} \leq \epsilon$ for each $J_{k}$.\vspace{5mm}

\noindent Now we will make use of a double-layered I-decomposition utilized in \cite{CR}. Subdivide $[0, T + \delta]$ into subintervals $J_{k}$, each $J_{k}$ is the union of $N^{1-}$ subintervals $J_{k,m}$ with $\| u_{\lambda} \|_{L_{t,x}^{4}(J_{k,m} \times \mathbf{R}^{3})} \leq \epsilon$ on each such subinterval. We will refer to the intervals $J_{k}$ as the big intervals, and the subintervals $J_{k,m}$ as the little intervals.\vspace{5mm}

\noindent Take the first big interval $J_{k}$. Crudely, by $(\ref{3.2})$, $E(Iu(t)) \leq 1$ on this big interval. Subdivide $J_{k} = \cup_{j = 0}^{N^{1-}} J_{k,m}$. Let $J_{k,m} = [a_{m}, b_{m}]$, $a_{0} = 0$, $a_{m + 1} = b_{m}$. The solution on $J_{k,m}$ will be written in the form

\begin{equation}\label{5.5}
e^{i(t - a_{m}) \Delta} u(a_{m}) + u_{j}^{nl}(t) = e^{it \Delta} u_{0} + \sum_{j = 1}^{m} e^{i(t - a_{j}) \Delta} u_{j - 1}^{nl}(a_{j}) + u_{m}^{nl}(t).
\end{equation}

\begin{equation}\label{5.6}
\sup_{t_{1}, t_{2} \in J_{k}} |E(Iu(t_{1})) - E(Iu(t_{2}))| \lesssim \frac{N^{1-}}{N^{2-}} + \frac{1}{N^{1-}} \| P_{> cN} \nabla Iu \|_{L_{t}^{2} L_{x}^{6}(J_{k} \times \mathbf{R}^{3})}^{2}.
\end{equation}

\noindent Now, by $(\ref{5.5})$,

\begin{equation}\label{5.7}
\aligned
\| P_{> cN} \nabla Iu \|_{L_{t}^{2} L_{x}^{6}(J \times \mathbf{R}^{3})} &\leq \| P_{> cN} \nabla Iu_{0} \|_{L_{x}^{2}(\mathbf{R}^{3})} + \sum_{m = 1}^{N^{1-}} \| \nabla P_{>cN} Iu_{m}^{nl}(a_{m}) \|_{L_{x}^{2}(\mathbf{R}^{3})} \\
&+ (\sum_{m = 0}^{N^{1-}} \| P_{> cN} \nabla Iu_{m}^{nl} \|_{L_{t}^{2} L_{x}^{6}(J_{k,m} \times \mathbf{R}^{3})}^{2})^{1/2}
\endaligned
\end{equation}

$$\| \nabla Iu_{0} \|_{L_{x}^{2}(\mathbf{R}^{3})} \lesssim 1,$$ which takes care of the first term. By $(\ref{4.7})$, $$\sum_{m = 1}^{N^{1-}} \| \nabla Iu_{m}^{nl}(a_{m}) \|_{L_{x}^{2}(\mathbf{R}^{3})} \lesssim \frac{N^{1-}}{N^{1-}} = 1,$$

\noindent which takes care of the second term. Finally,

$$(\sum_{m = 0}^{N^{1-}} \| P_{> cN} \nabla Iu_{m}^{nl} \|_{L_{t}^{2} L_{x}^{6}(J_{k,m} \times \mathbf{R}^{3})}^{2})^{1/2} \lesssim (\frac{N^{1-}}{N^{1-}})^{1/2} \lesssim 1.$$

\noindent In particular, this proves

\begin{equation}\label{5.8}
\sup_{t_{1}, t_{2} \in J_{k}} |E(Iu(t_{1})) - E(Iu(t_{2}))| \lesssim \frac{1}{N^{1-}}.
\end{equation}

\noindent When $s > 5/7$, $$CN^{\frac{3(1 - s)}{2s - 1}} << N^{2-},$$ so choosing $N$ sufficiently large proves

\begin{equation}\label{5.9}
\sup_{[0, T + \delta]} E(Iu_{\lambda}(t)) \leq \frac{9}{10}.
\end{equation}

\noindent This proves $W$ is both open and closed in $[0, \infty)$, so $W = [0, \infty)$.\vspace{5mm}

\noindent Finally, we prove scattering, following the argument in \cite{CKSTT2}. There is some $N$ such that

\begin{equation}\label{5.9.1}
E(Iu_{\lambda}(t)) \leq 1
\end{equation}

\noindent on $[0, \infty)$. By the interaction Morawetz estimates, $(\ref{5.4})$,

\begin{equation}\label{5.9.2}
\| u_{\lambda} \|_{L_{t,x}^{4}([0, \infty) \times \mathbf{R}^{3})} \leq C.
\end{equation}

\noindent Recall that by lemma $\ref{l2.1}$, if $\| u_{\lambda} \|_{L_{t,x}^{4}(J_{k, m} \times \mathbf{R}^{3})} \leq \epsilon$ and $E(Iu_{\lambda}(t)) \leq 1$ on $J_{k, m}$, then

\begin{equation}\label{5.10}
\| u \|_{L_{t}^{6} L_{x}^{9/2}(J_{k,m} \times \mathbf{R}^{3})} \lesssim (\epsilon^{2/3} + \frac{1}{N^{1/2}}).
\end{equation}

\noindent Let

\begin{equation}\label{5.11}
S_{s}(t) = \sup_{(p,q) \text{ admissible }} \| \langle \nabla \rangle^{s} u \|_{L_{t}^{p} L_{x}^{q}([0, t] \times \mathbf{R}^{3})}.
\end{equation}

$$S_{s}(t) \lesssim \| \langle \nabla \rangle^{s} u_{0} \|_{L^{2}(\mathbf{R}^{3})} + \| \langle \nabla \rangle^{s} u \|_{L_{t}^{2} L_{x}^{6}(J \times \mathbf{R}^{3})} \| u \|_{L_{t}^{6} L_{x}^{9/2}(J \times \mathbf{R}^{3})}^{2}$$

$$ \lesssim \| \langle \nabla \rangle^{s} u_{0} \|_{L^{2}(\mathbf{R}^{3})} + S_{s}(t) (\epsilon^{4/3} + \frac{1}{N}).$$

\noindent So for $\epsilon > 0$ sufficiently small and $N$ sufficiently large, this proves $S_{s}(t)$ is bounded on the first subinterval. Iterating over a finite number of subintervals proves $S_{s}(t) \leq C < \infty$ for $t \in [0, \infty)$. In particular, this proves

\begin{equation}\label{5.12}
\| u \|_{H^{s}(\mathbf{R}^{3})} \leq C(\| u_{0} \|_{H^{s}(\mathbf{R}^{3})}).
\end{equation}

\noindent Now set

\begin{equation}\label{5.13}
u_{+} = u_{0} + \int_{0}^{\infty} e^{-i \tau \Delta} |u(\tau)|^{2} u(\tau) d\tau.
\end{equation}

\begin{equation}\label{5.14}
\aligned
&\| \langle \nabla \rangle^{s} (e^{it \Delta} u_{+} - u(t,x)) \|_{L_{x}^{2}(\mathbf{R}^{3})} = \| \int_{t}^{\infty} \langle \nabla \rangle^{s} e^{-i \tau \Delta} |u(\tau)|^{2} u(\tau) d\tau \|_{L_{x}^{2}(\mathbf{R}^{3})} \\
&\lesssim \| \langle \nabla \rangle^{s} u \|_{L_{t,x}^{10/3}([T, \infty) \times \mathbf{R}^{3})} \| u \|_{L_{t,x}^{5}([T, \infty) \times \mathbf{R}^{3})}^{2}.
\endaligned
\end{equation}

\noindent As $T \rightarrow \infty$, $\| u \|_{L_{t,x}^{4}([T, \infty) \times \mathbf{R}^{3})} \rightarrow 0$, on the other hand,

\begin{equation}\label{5.15}
\| u \|_{L_{t,x}^{6}([0, \infty) \times \mathbf{R}^{3})} \lesssim \| \langle \nabla \rangle^{2/3} u \|_{L_{t}^{6} L_{x}^{18/7}([0, \infty) \times \mathbf{R}^{3})} \lesssim S_{2/3}(t) < \infty,
\end{equation}

\noindent by $(\ref{5.11})$. Interpolating proves $\| u \|_{L_{t,x}^{5}([T, \infty) \times \mathbf{R}^{3})} \rightarrow 0$ as $T \rightarrow \infty$. By Duhamel's principle

\begin{equation}\label{5.16}
\aligned
 \| \int_{0}^{\infty} &\langle \nabla \rangle^{s} e^{-i \tau \Delta} |u(\tau)|^{2} u(\tau) d\tau \|_{H^{s}(\mathbf{R}^{3})}	\lesssim	\| u \|_{L_{t,x}^{5}([0, \infty) \times \mathbf{R}^{3})}^{2} \| \langle \nabla \rangle^{s} u \|_{L_{t}^{5} L_{x}^{30/11}([0, \infty) \times \mathbf{R}^{3})}	\\
&\lesssim (\sup_{t \in [0, \infty)} S_{s}(t)) \| u \|_{L_{t,x}^{5}([0, \infty) \times \mathbf{R}^{3})}^{2} < \infty.
\endaligned
\end{equation}

\noindent Also,

\begin{equation}\label{5.17}
 \| \int_{T}^{\infty} e^{-i\tau \Delta} |u(\tau)|^{2} u(\tau) d\tau \|_{L_{x}^{2}(\mathbf{R}^{3})}	\lesssim	\| u \|_{L_{t,x}^{5}([T, \infty) \times \mathbf{R}^{3})}^{2} \| \langle \nabla \rangle^{s} u \|_{L_{t}^{5} L_{x}^{30/11}([T, \infty) \times \mathbf{R}^{3})}	\rightarrow 0
\end{equation}

 \noindent as $T \rightarrow \infty$. This completes the proof of theorem \ref{t1.1}. $\Box$

\newpage
\nocite*
\bibliographystyle{plain}
\bibliography{cubicn=3}
\end{document}